\documentclass{amsart}
\usepackage{amssymb}
\usepackage{amsfonts}
\usepackage{latexsym}

\newtheorem{theorem}{Theorem}[section]
\newtheorem{lemma}[theorem]{Lemma}
\newtheorem{proposition}[theorem]{Proposition}
\newtheorem{corollary}[theorem]{Corollary}
\theoremstyle{definition}

\newtheorem{remark}[theorem]{Remark}

\newcommand{\Tr}{\text{Tr}}

\newcommand{\End}{\text{End}}

\newcommand{\Rep}{\text{Rep}}

\newcommand{\ben}{\begin{enumerate}}
\newcommand{\een}{\end{enumerate}}

\begin{document}

\title{On Hopf Algebras of Dimension $pq$}
\author{Pavel Etingof}
\address{Department of Mathematics, Massachusetts Institute of Technology,
Cambridge, MA 02139, USA} \email{etingof@math.mit.edu}

\author{Shlomo Gelaki}
\address{Department of Mathematics, Technion-Israel Institute of
Technology, Haifa 32000, Israel}
\email{gelaki@math.technion.ac.il}

\maketitle

\section{Main theorem}

The classification of Hopf algebras of dimension $pq$, where $p$,
$q$ are distinct prime numbers, over an algebraically closed field
$k$ of characteristic $0$ is still unknown. In the previous works
\cite{eg}, \cite{gw} the problem was solved in the semisimple
case; namely, it was shown that any semisimple Hopf algebra of
dimension $pq$ is trivial (i.e. isomorphic to either a group
algebra or to the dual of a group algebra). In the work \cite{n},
Ng completed the classification of Hopf algebras of dimension
$p^2$, which started in \cite{as}, \cite{m}. In addition to that,
the problem was also settled in some specific low dimensions. More
precisely, Williams did dimensions 6 and 10 \cite{w},
Andruskiewitsch and Natale did dimensions 15, 21 and 35 \cite{an},
and Beattie and Dascalescu did dimensions 14, 55 and 77 \cite{bd}.

In this paper we contribute to the classification of Hopf algebras
of dimension $pq$ by proving the following main theorem. (We shall
work over the field of complex numbers, for convenience.)

\begin{theorem} Let $p<q$ be odd primes with $q\le 2p+1$.
Then any complex Hopf algebra of dimension $pq$ is semisimple and
hence trivial.
\end{theorem}

In particular, this theorem covers all the odd dimensions
listed above.

\section{Proof of the main theorem}

\begin{proposition}\label{estimate}
Let $H$ be a finite dimensional Hopf algebra.
Let $P$ denote the projective cover
of the trivial representation of $H$, and $N$ be the
number of irreducible representations of $H$.
Then $N\dim(P)\le \dim(H)$.
\end{proposition}

\begin{proof}
$H=\bigoplus_{V\in {\rm Irr}(H)}\dim(V)P(V)$, where $P(V)$ is the
projective cover of $V$. Now, the module $P(V)\otimes {}^*V$ is
projective (cf. \cite{kl}, p.441, Corollary 2 or \cite{gms},
Proposition 2.1) and ${\rm Hom}(P(V)\otimes {}^*V,\Bbb C)= {\rm
Hom}(P(V), V)\ne 0$. Hence $P(V)\otimes {}^*V$ contains $P$, and
its dimension is at least $\dim(P)$. So $\dim(V)\dim (P(V))\ge
\dim(P)$. Adding over all $V$, we get $N\dim(P)\le \dim(H)$, as
desired.
\end{proof}

\begin{lemma}\label{PC} If $P=\Bbb C$ then
$H$ is semisimple.
\end{lemma}

\begin{proof} If $P=\Bbb C$ then $\Bbb C$ is
projective, and hence any representation $V=V\otimes \Bbb C$ of
$H$ is projective (\cite{kl}, p.441, Corollary 2; \cite{gms},
Proposition 2.1). This implies that $\Rep (H)$ is a semisimple
category, so $H$ is semisimple.
\end{proof}

Let $H$ be a {\bf non-semisimple} complex Hopf algebra of
dimension $pq$, $p<q$ are primes.

\begin{lemma}\label{1d}
Let $D\subset \Rep(H)$ be the subcategory of representations whose
composition factors are 1-dimensional. Then $D$ is semisimple.
\end{lemma}

\begin{proof} Assume that $D$ is not semisimple.
Then $D$ contains a 2-dimensional representation $V$ of $H$ which
is a nontrivial extension of a character $\chi_1$ by a character
$\chi_2$. This means, $V$ has a basis $v_1,v_2$ such that for any
$a\in H$, one has $av_1=\chi_1(a)v_1+f(a)v_2$,
$av_2=\chi_2(a)v_2$. Since $(ab)v_1=a(bv_1)$, the element $f\in
H^*$ satisfies the equality $f(ab)=f(a)\chi_1(b)+\chi_2(a)f(b)$.
Thus, $\Delta(f)= f\otimes \chi_1+\chi_2\otimes f$, so $f$ is
skew-primitive. Since $V$ is a nonsplit extension, $f$ is
nontrivial (i.e. not a multiple of $\chi_1-\chi_2$). Thus, $H^*$
contains nontrivial skew-primitive elements, which contradicts
Proposition 1.8 in [AN]. Thus $D$ is semisimple.
\end{proof}

Suppose from now on that $p$ and $q$ are {\bf odd}.

\begin{theorem}\label{ng} \cite{n} (i) The order of $S^4$ is
$p$.

(ii) The distinguished grouplike element either in $H$ or in
$H^*$ is nontrivial and has order $p$.
\end{theorem}

So we will assume that $H^*$ has distinguished grouplike element
$\chi$ of order $p$. Then $\chi$ is a 1-dimensional representation
of $H$, $\chi^p=1$.

\begin{proposition}\label{stab}
Let $V$ be an irreducible $H$-module, which is stable under
tensoring on the left with $\chi$. Then $\dim(V)$ is divisible by
$p$.
\end{proposition}

\begin{proof} As pointed out in the proof
of Lemma 2.1 in [AN], $(\End V)^*\subset H^*$ is a Hopf module for
$\Bbb (\mathbb{C}[G(H^*)],H^*)$. Thus the result follows by the
Nichols-Zoeller theorem \cite{nz}.
\end{proof}

\begin{proposition}\label{q2p}
Suppose that $V$ is an irreducible $H$-module, which is stable
under left multiplication by $\chi$. Then $q> 2p+1$.
\end{proposition}

\begin{proof}
Assume that $q\le 2p+1$. By Proposition \ref{stab}, the dimension
of $V$ is divisible by $p$. We claim that $V$ is the unique
irreducible $H$-module of dimension $\ge p$, and hence the unique
irreducible module stable under left multiplication by $\chi$; in
particular $V^*={}^*V=V$. Indeed, assume that $W$ is another
irreducible $H$-module of dimension $\ge p$. Then $$ pq=\dim(H)\ge
\dim(V)^2+\dim(W)^2+\sum_{j=1}^p \dim(\chi^j)^2=2p^2+p. $$ Since
$q\le 2p+1$, this inequality is an equality, and hence $H$ is
semisimple, a contradiction.

Assume that $Y,Z$ are other (hence $\chi$-unstable) simple modules
such that $Y\otimes Z$ contains $V$ as a constituent of the
Jordan-H\"older series. Then
$\dim(Y)>1$, $\dim(Z)>1$, and either $Y$ or $Z$ (say, $Y$) has
dimension $> \sqrt{p}$, so $$ pq=\dim (H)\ge \sum_{j=1}^p
\dim(\chi^j \otimes Y)^2+(\dim (V))^2+\sum_{j=1}^p\dim(\chi^j)^2>
2p^2+p, $$ i.e. $q>2p+1$, a contradiction. Thus, $Y,Z$ do not
exist.

Let $C$ be the category of $H$-modules which do not contain $V$ as
a constituent. As we have seen in the previous paragraph, it is a
rigid tensor category, so $C=\Rep H'$, where $H'$ is a quotient of
$H$ by a nontrivial Hopf ideal. Since $\chi^j$, $1\le j\le p$, are
representations of $H'$, by the Nichols-Zoeller theorem \cite{nz},
the dimension of $H'$ is a divisor of $pq$ which is divisible by
$p$. Hence the dimension of $H'$ is $p$, and we have $H'=\Bbb
C[\Bbb Z_p]$. In other words, $C$ is semisimple, with simple
objects $\chi^j$, $1\le j\le p$.

Let $P(V)$ be the projective cover of $V$. If $P(V)\ne V$ then the
socle of $P(V)$ is $\chi^{-1}\otimes V=V$ (by [EO], Lemma 2.10).
Thus, $\dim (P(V))\ge 2\dim (V)$, and $$ pq=\dim(H)\ge \dim(V)\dim
(P(V))+\sum_{j=1}^p\dim(\chi^j)\dim (P(\chi^j))> 2p^2+p $$ (the
$>$ sign is due to the fact that the projective cover of $\Bbb C$
cannot equal $\Bbb C$, by Lemma \ref{PC}).
So, we have a contradiction. Thus $P(V)=V$.

This means that $P(\Bbb C)$ must contain $V$ as a constituent
(otherwise by Lemma \ref{1d}, one has $P(\Bbb C)=\Bbb C$, hence
by Lemma \ref{PC} $H$
 is semisimple, which is a contradiction).
Hence $P(\chi^j)$ contains $V$ as a constituent for all $j$.
Thus the
dimension of $P(\chi^j)$ is at least $p+2$ (it involves $\chi^j$,
$V$, and the socle $\chi^{j-1}$). So $$ pq=\dim (H)\ge
\sum_{j=1}^p \dim (P(\chi^j))+(\dim (V))^2\ge p(p+2)+p^2= p(2p+2),
$$ a contradiction.
\end{proof}

Assume from now on that ${\bf q\le 2p+1}$.

\begin{lemma} \label{dimp}
One has, $\dim(P)<p$.
\end{lemma}

\begin{proof}
By Proposition \ref{q2p}, all simple $H$-modules are
$\chi$-unstable. By Lemma \ref{1d}, not all simple $H$-modules are
1-dimensional. Thus there are at least $2$ orbits of simple
$H$-modules under left multiplication by $\chi$. Consider two
cases.

1. The number of orbits is $\ge 3$. Then by
Proposition \ref{estimate}, $\dim (P)\le q/3\le \frac{2p+1}{3}<
p$.

2. The number of orbits is $2$. Then by Proposition
\ref{estimate}, $\dim(P)\le q/2$, so $\dim(P)\le p$. So it remains
to show that $\dim(P)\ne p$. Suppose $\dim(P)=p$. Let $V$ be an
element of the unique orbit of non-1-dimensional representations,
and suppose $\dim(V)=d$. Then $P(V)\otimes {}^*V$ contains $P$, so
$\dim (P(V))\ge p/d$. Clearly, we cannot have $P=P(V)\otimes
{}^*V$, as the dimension of $P$ is a prime. Thus, $P(V)\otimes
{}^*V=P\oplus Q$, where $Q$ is a nontrivial projective module.
Thus, $\dim(Q)\ge p/d$ (as $Q$ contains $P(\chi^j)=\chi^j\otimes
P$, or it contains $P(\chi^j\otimes V)=\chi^j\otimes P(V)$,
because they are the only indecomposable projective $H$-modules).
So $\dim (P(V))\dim(V)\ge p+p/d=p(d+1)/d$. This implies (like in
the proof of Proposition \ref{estimate}) that
\begin{eqnarray*} \lefteqn{pq\ge \sum_{j=1}^p \dim
(P(\chi^j\otimes V))\dim (\chi^j\otimes V)+ \sum_{j=1}^p \dim
(P(\chi^j))}\\ & & \ge p^2(d+1)/d+p^2=p^2\frac{2d+1}{d}>p(2p+1)
\end{eqnarray*}
(since clearly $d<p$: otherwise $pq=\dim(H)\ge p\dim(V)^2=p^3$,
which is absurd). This is a contradiction, so the lemma is proved.
\end{proof}

\begin{lemma} \label{conjchi}
For any simple $H$-module $V$, one has $\chi^{-1}\otimes
V\otimes \chi=V$.
\end{lemma}

\begin{proof}
Let $Y,Z$ be irreducible representations of $H$ which are stable
under conjugation by $\chi$, and suppose that $Y\otimes Z$
contains as a constituent an irreducible representation $X$ which
is not stable under conjugation by $\chi$. Then $Y\otimes Z$ must
also contain as constituents all the translates $\chi^i\otimes
X\otimes \chi^{-i}$, which are pairwise non-isomorphic. To be
definite, assume that $\dim(Y)\ge \dim(Z)$. Since $X$ cannot be
1-dimensional, we find that $\dim(Y)\dim(Z)\ge 2p$, so $\dim(Y)^2>
2p$. This means that $$ pq=\dim (H)\ge \sum_{j=1}^p
(\dim(\chi^j\otimes Y)^2+\dim(\chi^j)^2)> 2p^2+p, $$ a
contradiction. Thus $Y,Z$ do not exist.

So let $E$ be the category of representations of $H$ whose all
simple constituents are stable under conjugation by $\chi$. As we
have just shown, it is a rigid tensor category. So $E=\Rep(H')$,
where $H'$ is a quotient of $H$.

Assume $H'\ne H$.
Then $H'=\Bbb C[\Bbb Z_p]$.
This means that all simple modules
over $H$ which are not 1-dimensional fail
to be invariant under conjugation by $\chi$.

Choose a non-1-dimensional simple constituent $V$ in $P$. Such
exists, otherwise by Lemma \ref{1d}, $\Bbb C$ is projective and
$H$ is semisimple, which is a contradiction. Then all translates
$\chi^i\otimes V\otimes \chi^{-i}$ (which are pairwise
non-isomorphic) are contained in $P$ as constituents. Thus
$\dim(P)\ge 2p$ (since $\dim(V)\ge 2$).
But we showed before that $\dim(P)<p$.
Contradiction. Thus, $H=H'$, as desired.
\end{proof}

\begin{corollary}\label{c1}
For any simple $H$-module $V$, one has
$V=V^{****}$.
\end{corollary}

\begin{proof}
Follows from Lemma \ref{conjchi} and the Radford $S^4$ formula
\cite{r}.
\end{proof}

\begin{corollary}\label{c2}
There exists an odd dimensional
indecomposable projective module
$Q$ over $H$ such that $Q=Q^{**}$.
\end{corollary}

\begin{proof}
We have $H=\bigoplus \dim(V)P(V)$ and $H^{**}=H$. By Corollary
\ref{c1}, $**$ is an involution on the set of indecomposable
projectives $P(V)$. Since the dimension of $H$ is odd, the sum
$\bigoplus_{V=V^{**}}\dim(V)P(V)$ is also odd dimensional. Thus
one of the summands is odd dimensional, as desired.
\end{proof}

\begin{lemma}\label{oddim}
Let $Q$ be an odd dimensional indecomposable
projective module over $H$ such that $Q=Q^{**}$. Then $\dim(Q)\ge p$.
\end{lemma}

\begin{proof} Let $a: Q\to Q^{**}$ be an isomorphism.
In other words, $a$ is a linear operator $Q\to Q$
such that $ax=S^2(x)a$, $x\in H$.
We claim that if $\lambda,\mu$ are eigenvalues
of $a$ then $\lambda/\mu$ is
a root of unity of degree $2p$, i.e. a power of $z=e^{\pi i/p}$.

Indeed, assume that this is not the case.
Let $Q_1$ be the sum of generalized eigenspaces of $a$ on $Q$
with eigenvalues $\lambda z^m$, $m=0,...,2p-1$, and $Q_2$ the sum
of generalized
eigenspaces of $a$ with other eigenvalues. Then $Q_1,Q_2\ne 0$, and
$Q_1\oplus Q_2=Q$. We claim that $Q_1$ and $Q_2$ are
$H$-submodules of $Q$ (this is a contradiction, since $Q$ is
indecomposable). Indeed, let $x$ be an eigenvector of $S^2$
in $H$ with eigenvalue $z^m$, and $v\in Q_i$ for some $i$ be
a generalized eigenvector of $a$ with eigenvalue $\beta$:
 $(a-\beta)^Nv=0$. Since $ax=S^2(x)a=z^mxa$, we have
$0=x(a-\beta)^Nv=(z^{-m}a-\beta)^Nxv$.
This shows that $xv$ is a generalized eigenvector of
$a$ with eigenvalue $\beta z^m$, i.e. $xv\in Q_i$ for the same
$i$. But by Theorem \ref{ng}, elements $x$ as above span $H$.
We are done.

Let us normalize $a$ in such a way that one of its eigenvalues is
$1$. Then, as we have shown, other eigenvalues will be $z^j$,
$j=0,...,2p-1$.

By Lorenz' theorem \cite{l}, (see also, \cite{eo}, proof of
Theorem 2.16), $\Tr(a)=0$. Thus, if $m_j$ is the multiplicity of
the eigenvalue $z^j$, then $$ \sum_{j=0}^{2p-1} m_jz^j=0,\;
\sum_{j=0}^{2p-1} m_j=\dim(Q). $$ Since $z^{p+j}=-z^j$, we can
rewrite the last equation as $$
\sum_{j=0}^{p-1}(m_{2j}-m_{2j+p})z^{2j}=0. $$ This means that
$m_{2j}-m_{2j+p}=b$ for all $j$. Since $\sum_j m_j$ is odd, $b\ne
0$. Thus, ${\rm max}(m_{2j},m_{2j+p})\ge 1$, and hence $\dim(Q)\ge
p$, as desired.
\end{proof}

\begin{remark}
The proof of Lemma \ref{oddim} is similar to \cite{n}, Section 5.
\end{remark}

\begin{lemma}\label{Podd}
The module $P$ is odd dimensional.
\end{lemma}

\begin{proof} By Corollary \ref{c2} and Lemma \ref{oddim}, some
indecomposable projective module $Q=P(V)$ over $H$ has dimension
$\ge p$. If $V$ is not 1-dimensional, then $$ pq=\dim (H)\ge
\sum_{j=1}^p \dim (P(\chi^j\otimes V))\dim (\chi^j\otimes V)+
\sum_{j=1}^p \dim (P(\chi^j))\dim(\chi^j)> 2p^2+p$$ (since
$\dim(V)\ge 2$ and $P(\Bbb C)\ne \Bbb C$ by Lemma \ref{PC}), a
contradiction. So $V$ is 1-dimensional, hence we may assume that
$V=\Bbb C$, as desired.
\end{proof}

Now we observe that Lemma \ref{Podd} and Lemma \ref{oddim} imply
that $\dim P\ge p$, which contradicts our previous conclusion
(Lemma \ref{dimp}) that $\dim (P)<p$. Thus we get an ultimate
contradiction, which means that nonsemisimple Hopf algebras of
dimension $pq$ with $p<q\le 2p+1$, where $p$, $q$ are primes, do
not exist. The main theorem is proved.

\section{Acknowledgments}
The authors thank V. Ostrik for useful discussions.
P.E was partially supported by the NSF grant DMS-9988796. P.E
partially conducted his research for the Clay Mathematics
Institute as a Clay Mathematics Institute Prize Fellow. S.G's
research was supported by Technion V.P.R. Fund - Dent Charitable
Trust- Non Military Research Fund, and by THE ISRAEL SCIENCE
FOUNDATION (grant No. 70/02-1).


\begin{thebibliography}{AEG}

\bibitem[AN]{an} N. Andruskiewitsch and S. Natale, Counting arguments
for Hopf algebras of low dimension, {\em Tsukuba J. Math.} {\bf
25} (2001), no. 1, 187--201.

\bibitem[AS]{as} N. Andruskiewitsch and H-J. Schneider, Hopf algebras
of order $p\sp 2$ and braided Hopf algebras of order $p$ {\em J.
Algebra} {\bf 199} (1998), no. 2, 430--454.

\bibitem[BD]{bd} M. Beattie and S. Dascalescu, Hopf algebras of dimension
14, {\em preprint}, math.QA/0205243.

\bibitem[EG]{eg} P. Etingof and S. Gelaki, Semisimple Hopf algebras of
dimension $pq$ are trivial, {\em Journal of Algebra} {\bf 210}
(1998), 664--669.

\bibitem[EO]{eo} P. Etingof and V. Ostrik, Finite tensor
categories, {\em preprint}, math.QA/0301027.

\bibitem[GMS]{gms} E. L. Green, E. N. Marcos and Solberg, Representations
and almost split sequences for Hopf algebras, Representation
theory of algebras (Cocoyoc, 1994), 237--245, CMS Conf. Proc.,
{\bf 18}, Amer. Math. Soc., Providence, RI, 1996.

\bibitem[GW]{gw} S. Gelaki and S. Westreich, On semisimple Hopf
algebras of dimension $pq,$ {\em Proceedings of the AMS}, {\bf
128} (2000), no.1, 39--47.

\bibitem[KL]{kl} D. Kazhdan and G. Lusztig,  Tensor structures arising
from affine Lie algebras. IV. {\em J. Amer. Math. Soc.} {\bf 7}
(1994), no. 2, 383--453.

\bibitem[L]{l} M. Lorenz, Representations of finite-dimensional Hopf
algebras, {\em J. Alg.} {\bf 188} (1997), no. 2, 476--505.

\bibitem[M]{m} A. Masuoka, The $p\sp n$ theorem for semisimple Hopf algebras,
{\em Proc. Amer. Math. Soc.} {\bf 124} (1996), no. 3, 735--737.

\bibitem[N]{n} S-H. Ng, Non-semisimple Hopf algebras of dimension
$p\sp 2$ {\em J. Algebra} {\bf 255} (2002), no. 1, 182--197.

\bibitem [NZ]{nz} W.D. Nichols and M.B. Zoeller, A Hopf algebra freeness
theorem, {\em American Journal of Mathematics} {\bf 111} (1989),
381--385.

\bibitem[R]{r} D.E. Radford, The order of the antipode of a
finite-dimensional Hopf algebra is finite, {\em Amer. J. Math.}
{\bf 98} (1976), 333--355.

\bibitem[W]{w} R. Williams, Finite Dimensional Hopf algebras, {\em
Ph.D thesis, Florida State University}, 1998.

\end{thebibliography}
\end{document}